\documentclass[reqno,11pt]{amsart}

\usepackage{amsthm, amsmath, amssymb, stmaryrd}
\usepackage[utf8]{inputenc}
\usepackage[T1]{fontenc}

\usepackage[hyphenbreaks]{breakurl}
\usepackage[hyphens]{url}

\usepackage[shortlabels]{enumitem}
\usepackage[hidelinks]{hyperref}
\usepackage{microtype}

\usepackage{bm}
\usepackage[margin=1in]{geometry}

\usepackage[textsize=scriptsize,backgroundcolor=orange!5]{todonotes}

\usepackage[noabbrev,capitalize,sort]{cleveref}
\crefname{equation}{}{}
\numberwithin{equation}{section}

\newtheorem{theorem}{Theorem}[section]
\newtheorem{proposition}[theorem]{Proposition}
\newtheorem{lemma}[theorem]{Lemma}
\newtheorem{claim}[theorem]{Claim}
\newtheorem{corollary}[theorem]{Corollary}

\newtheorem*{theorem11}{Theorem 1.1}
\newtheorem*{theorem12}{Theorem 1.2}

\theoremstyle{definition}
\newtheorem{definition}[theorem]{Definition}

\newtheorem{question}[theorem]{Question}
\newtheorem{example}[theorem]{Example}

\theoremstyle{remark}
\newtheorem*{remark}{Remark}

\newcommand{\abs}[1]{\left\lvert#1\right\rvert}

\DeclareMathOperator{\Span}{span}
\DeclareMathOperator{\tr}{tr}

\newcommand{\RR}{\mathbb{R}}

\newcommand{\ZZ}{\mathbb{Z}}

\newcommand{\cF}{\mathcal F}

\newcommand{\cA}{\mathcal A}

\newcommand{\bS}{\mathbf S}
\newcommand{\bT}{\mathbf T}

\newcommand*{\eqdef}{\stackrel{\mbox{\normalfont\tiny def}}{=}}   
\newcommand*{\PP}{\mathbb{P}}  

\newcommand\coin{
\mathrel{\ooalign{\hss$\bigcirc$\hss\cr\kern0.7ex\hbox{\scalebox{0.8}{$\$$}}}}\,}

\makeatletter
\newcommand\thankssymb[1]{\textsuperscript{\@fnsymbol{#1}}}
\makeatother

\author[Ting-Wei Chao]{Ting-Wei Chao\thankssymb{1}}
\author[Hung-Hsun Hans Yu]{Hung-Hsun Hans Yu\thankssymb{2}}

\thanks{\thankssymb{1}Department of Mathematical Sciences, Carnegie Mellon University, Pittsburgh, PA 15213, USA\@. Supported in part
by U.S. taxpayers through NSF grant DMS-2154063 and NSF CAREER grant DMS-1555149. Email: {\tt tchao2@andrew.cmu.edu}}

\thanks{\thankssymb{2}Department of Mathematics, Princeton University, Princeton, NJ 08544\@.  Email: {\tt hansonyu@princeton.edu}}

\title{Kruskal--Katona-Type Problems via the Entropy Method}

\begin{document}

\maketitle

\begin{abstract}
	In this paper, we investigate several extremal combinatorics problems that ask for the maximum number of copies of a fixed subgraph given the number of edges. We call problems of this type Kruskal--Katona-type problems.
 Most of the problems that will be discussed in this paper are related to the joints problem. There are two main results in this paper. First, we prove that, in a $3$-edge-colored graph with $R$ red, $G$ green, $B$ blue edges, the number of rainbow triangles is at most $\sqrt{2RGB}$, which is sharp. Second, we give a generalization of the Kruskal--Katona theorem that implies many other previous generalizations. Both arguments use the entropy method, and the main innovation lies in a more clever argument that improves bounds given by Shearer's inequality.
\end{abstract}

\section{Introduction}
We begin with a brief introduction to the joints problem as it motivates most graph-theoretic problems in this paper.
The joints problem asks to find the maximum number of joints formed by $N$ lines in $\RR^d$, where a joint is a point that lies on $d$ lines whose directions are linearly independent.
First studied by Chazelle et al. \cite{CEGPSSS92}, it has drawn attention not only because it is an interesting incidence geometry problem on its own, but also because there is a connection between the joints problem and the Kakeya problem in harmonic analysis as observed by Wolff \cite{Wol99}.
For a more comprehensive review of works related to the joints problem and other variants, we refer the readers to \cite{TYZ22}, where several generalizations are also discussed.
We will focus on the connection between the joints problem and some extremal problems in graph theory.

In the original paper \cite{CEGPSSS92} where the joints problem was proposed, a construction of $N$ lines with many joints was given as follows.
Choose $k$ hyperplanes in $\RR^d$ in general position, and consider the $(d-1)$-wise intersections of the $k$ hyperplanes, giving $N=\binom{k}{d-1}$ lines.
It is easy to verify that the $d$-wise intersections of the hyperplanes are joints formed by the $N$ lines, giving $\binom{k}{d}$ joints.
Guth \cite[Section 2.5]{Guth-book} conjectured that this is optimal, and this was verified asymptotically in a work of Zhao and the second author \cite{YZ23}.
We have recently verified this conjecture exactly in \cite{CY23}.

The specific construction above suggests a more general family of constructions for the joints problem using hyperplanes in general position.
Indeed, we may consider a $(d-1)$-uniform hypergraph $H$ where each vertex represents a hyperplane, and each hyperedge represents the line $\ell$ that is the intersection of the hyperplanes represented by the verices in the hyperedge.
A joint then corresponds to a $d$-set of the vertices all of whose $(d-1)$-subsets are present in $H$.
This type of configuration was already considered in \cite{YZ23}, and as in that paper, we call this type of configuration \emph{generically induced}.

When restricted to the generically induced configurations, the joints problem simplifies as follows.
First note that finding the maximum number of joints given the number of lines is the same as finding the minimum number of lines given the number of joints.
Suppose that the configuration is generically induced by a hypergraph $H$, and let $\cA$ be the family of $d$-subsets of $V(H)$ corresponding to the joints.
It is then clear that $\partial\cA\subseteq E(H)$ (where $\partial \cA\eqdef\{S-\{s\}\mid S\in \cA, s\in S\}$ is the lower shadow), and having those suffices to form all joints in $\cA$.
Therefore the well-known Kruskal--Katona theorem tells us exactly the minimum number of lines given the number of joints in the generically induced case.

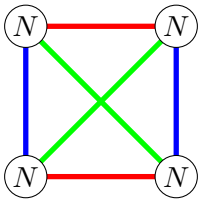
\begin{figure}[h]
    \centering
    
\begin{tikzpicture}
    \coordinate (A) at (0,0);
    \coordinate (B) at (2,0);
    \coordinate (C) at (2,2);
    \coordinate (D) at (0,2);
    
    \draw[red, line width=2pt] (A) -- (B);
    \draw[red, line width=2pt] (C) -- (D);
    \draw[blue, line width=2pt] (A) -- (D);
    \draw[blue, line width=2pt] (B) -- (C);
    \draw[green, line width=2pt] (A) -- (C);
    \draw[green, line width=2pt] (B) -- (D);
    
    \foreach \vertex/\label in {A/$N$, B/$N$, C/$N$, D/$N$}
        \filldraw[fill=white] (\vertex) circle (8pt) node {\label};
\end{tikzpicture}
    \caption{$2N^2$ edges of each color with $4N^3$ rainbow triangles.}
    \label{fig:rt}
\end{figure}

Many other variants of the joints problems have also been considered, and in those cases, not much is known about the exact answer even in the generically induced cases.
A representative case is the multijoints problem in $\RR^3$, where instead of a single set of lines, we now have three sets of lines of sizes $L_1,L_2$ and $L_3$, respectively.
A \emph{multijoint} is a joint that is formed by three lines coming from three different sets.
For a generically induced configuration corresponding to a graph $G$, we may think of $G$ as a $3$-edge-colored graph where the numbers of the edges of the three colors are $L_1,L_2$ and $L_3$, respectively.
In this language, a multijoint is a $3$-subset of the vertices that forms a rainbow triangle.
The best known upper bound for the general multijoints problem was proved in \cite{YZ23}, showing that the number of multijoints is at most $\sqrt{6L_1L_2L_3}$.
However, the best known construction, which appeared in \cite{YZ23}, only gives $4N^3$ multijoints with $2N^2$ lines of each color with a generically induced configuration.
In the language of rainbow triangles, the construction is obtained by blowing up a $K_4$ where opposite edges have the same color, see \cref{fig:rt}.
As mentioned in \cite{YZ23}, a preliminary flag algebra computation was done by Lidick\'{y} that showed an upper bound of $\sqrt{2N^3}$ when $L_1=L_2=L_3=N$ in the generically induced case.
One of our main results is to verify this for all the generically induced configurations using the entropy method.

\begin{theorem}\label{theorem:RainbowTriangle}
Let $G=(V,E)$ be a simple graph and let edges be colored in colors red, green, and blue, so that each edge receives exactly one color. Let $R,G,B$ denote the number of red, green, blue edges, respectively. Let $T$ be the number of rainbow triangles in the graph.
Then $T^2\leq 2RGB$. 
\end{theorem}

Balogh et al. \cite{BHLPVY17} tackled a similar problems with flag algebra, though there the number of vertices is given instead of the number of edges.

We shall comment that this type of problem had already been considered in much more generality even before its connection to the joints problem was noticed---extremal problems that ask for the maximum number of copies of a subgraph given the number of edges were addressed by Alon \cite{Alon81}, and the hypergraph analog was studied by Friedgut and Kahn \cite{FK98}.
As the Kruskal--Katona theorem seems to be the most known instance of this type of problem, we refer to those problems as \emph{Kruskal--Katona-type problems}.
Although the two aforementioned works have determined the extremal number up to a constant (depending on the subgraph one counts), less is known about the precise asymptotic behavior and the multicolored variant of the problem.

\cref{theorem:RainbowTriangle} is proved using the entropy method.
This should not be a surprise as the previous work of Friedgut and Kahn \cite{FK98} also used the entropy method.
In their work, the main tool was Shearer's inequality, which has proven to be useful in many applications of the entropy method.
In the case of the rainbow triangle problem, Shearer's inequality immediately gives the right growth rate of the maximum number of rainbow triangles, though it gives an upper bound that is off from the answer by a factor of two (we will discuss this formally in \cref{section:RainbowTriangle}).
Closing this gap requires us to not use Shearer's inequality as a black box, and in this case it requires some new idea.
We think this is interesting as there are not many known instances where effort is put to beat a bound given by Shearer's inequality with a more clever entropy argument.
A near-example that immediately comes to mind is the work of Chung et al. \cite{CGFS86} that used Shearer's inequality to bound the size of triangle-intersecting family of graphs, where the gap of a factor of two was later closed up by Ellis, Filmus and Friedgut \cite{EFF12}, though they used spectral methods instead.

The Kruskal--Katona theorem provides yet another setting where there is a quick way to get an upper bound using Shearer's inequality that is tight only up to an $o(1)$-factor.
At first sight, it seems difficult to prove a more accurate upper bound, as the Kruskal--Katona bound does not seem compatible with the entropy method.
The second topic of this paper is an entropy argument that gives a new proof of the so-called Lov\'asz's version of the Kruskal--Katona theorem.

\begin{theorem}[Lov\'asz's version of the Kruskal--Katona theorem]\label{theorem:Lovasz}
    Let $\cA$ be a family of $d$-subsets of $[n]$ and $\partial \cA$ be all the $(d-1)$-subsets that are contained in some sets in $\cA$. If $\abs{\cA}=\binom{t}{d}$ for some real number $t\geq d$, then $\abs{\partial\cA}\geq \binom{t}{d-1}$.
\end{theorem}

In fact, we prove a much more general result that vastly generalizes Lov\'asz's version of the Kruskal--Katona theorem.
As the statement is technical and hard to state, we postpone the statement until \cref{section:KruskalKatona}.
This generalization also implies some other known closely-related results, including the $q$-analog of \cref{theorem:Lovasz}.

\begin{theorem}[Chowdhury--Patk\'{o}s \cite{CP10}]\label{theorem:q-analog}
Let $\cA$ be a family of $d$-dimensional subspaces of $\mathbb{F}_q^n$ and $\partial \cA$ be all the $(d-1)$-dimensional subspaces that are contained in some subspace in $\cA$. If $t\geq d$ is a real number that satisfies
\[\abs{\cA}=\bigl[\!\begin{smallmatrix} t \\ d \end{smallmatrix}\!\bigr]_q\eqdef \frac{(q^t-1)(q^t-q)\ldots(q^t-q^{d-1})}{(q^d-1)(q^d-q)\ldots(q^d-q^{d-1})},\]
then 
\[\abs{\partial \cA}\geq \bigl[\!\begin{smallmatrix} t \\ d-1 \end{smallmatrix}\!\bigr]_q=\frac{(q^t-1)(q^t-q)\ldots(q^t-q^{d-2})}{(q^{d-1}-1)(q^{d-1}-q)\ldots(q^{d-1}-q^{d-2})}.\]

\end{theorem}

This paper is structured as follows.
In \cref{section:Prelim}, we introduce the definition and several properties of entropy that will be used throughout the paper.
In \cref{section:RainbowTriangle}, we give a proof of \cref{theorem:RainbowTriangle} and briefly discuss some partial results about the natural generalization of the theorem.
In \cref{section:KruskalKatona}, we give an entropy-based proof of \cref{theorem:Lovasz} and \cref{theorem:q-analog} by showing a common generalization.
Lastly, in \cref{section:OpenProblem}, we discuss some other extremal-graph-theoretic problems that are inspired by the joints problem.

\section{Preliminaries}\label{section:Prelim}
In this section, we list the properties of entropy that will be used in the later sections.
For a complete discussion and proofs for those, we refer the readers to \cite[Section 14.6]{AS00}.

\begin{definition}
    For any discrete random variable $X$ with finite support $S$, its entropy is defined as
    \[H(X)\eqdef \sum_{x\in S}-p_X(x)\log_2p_X(x)\]
    where we denote by $p_X(x)$ the probability $\PP(X=x)$.
    
    The entropy $H(X_1,\ldots, X_n)$ of several discrete random variables $X_1,\ldots, X_n$ is defined similarly.
    To be more specific, it is defined as the entropy of $\mathbf{X}=(X_1,\ldots, X_n)$.
\end{definition}

For example, a uniform random bit has entropy $1$.
More generally, a Bernoulli random variable with parameter $p$ has entropy $h(p)\eqdef -p\log_2p-(1-p)\log_2(1-p)$.
Informally speaking, the entropy measures the expected amount of information one gets by revealing the random variable, measured in bits.

The following proposition says that given the support of a random variable, its entropy is maximized by the uniform distribution.

\begin{proposition}\label{proposition:UnifBound}
    Let $X$ be a random variable with finite support $S$.
    Then $H(X)\leq \log_2\abs{S}$.
    The equality is attained when the law of $X$ is the uniform distribution on $S$.
\end{proposition}

There are several inequalities that we will use to bound entropies of the form $H(X_1,X_2,\ldots, X_n)$.
The first is subadditivity.

\begin{proposition}[Subadditivity]
    For any list of random variables $X_1,X_2,\ldots, X_n$ with finite support,
    \[H(X_1,X_2,\ldots,X_n)\leq \sum_{i=1}^{n}H(X_i).\]
    The equality holds if and only if all the random variables are mutually independent.
\end{proposition}

In fact, a much more general inequality also holds.
This is usually called Shearer's inequality.

\begin{theorem}[Shearer's inequality]\label{theorem:ShearerStrong}
    Let $k$ be a positive integer.
    Let $X_1,\ldots, X_n$ be a list of random variables with finite support, and let $I_1,\ldots,I_m$ be subsets of $[n]$ so that for each $i\in [n]$, $i\in I_j$ holds for at least $k$ different $j$'s. 
    Then
    \[kH(X_1,\ldots, X_n)\leq \sum_{j=1}^{m}H\left((X_i)_{i\in I_j}\right).\]
\end{theorem}

The following special case of Shearer's inequality would already be useful, and this is what we will apply in most of the cases in this paper.

\begin{theorem}[A special case of Shearer's inequality]\label{theorem:Shearer}
    Let $X_1,\ldots, X_n$ be a list of random variables with finite support.
    Then
    \[(n-1)H(X_1,\ldots, X_n)\leq \sum_{i=1}^{n}H(X_1,\ldots, \widehat{X_i},\ldots,X_n),\]
    where $\widehat{X_i}$ means omitting $X_i$ from the list.
\end{theorem}

The final concept we need is the conditional entropy.
Roughly speaking, it computes the expected entropy of one random variable when conditioned on the other random variable.

\begin{definition}
    Let $X,Y$ be two random variables with supports $S,T$, respectively, that are both finite.
    We still denote by $p_{X,Y}(x,y)$ and $p_Y(y)$ the probabilities $\PP(X=x, Y=y)$ and $\PP(Y=y)$, respectively.
    For each $y\in T$, let $X\mid Y=y$ be the conditional random variable $X$ given $Y=y$, and let $H(X\mid Y=y)$ be its entropy. Moreover, we say that $H(X\mid Y=y)=0$ if $y$ is not in the support of $Y$.
    With those notations, the conditional entropy of $X$ given $Y$ is defined as
    \[H(X\mid Y)\eqdef \sum_{(x,y)\in S\times T}-p_{X,Y}(x,y)\log_2\left(\frac{p_{X,Y}(x,y)}{p_Y(y)}\right)=\sum_{y\in T}-p_Y(y)H(X\mid Y=y).\]
\end{definition}

With conditional entropy introduced, we may now state the chain rule.

\begin{proposition}[Chain rule]
    For any two random variables $X,Y$ with finite support, we have $H(X,Y) = H(Y)+H(X\mid Y)$.
    More generally, for any random variables $X_1,\ldots, X_n$ all with finite support,
    \[H(X_1,\ldots, X_n) = \sum_{i=1}^{n}H(X_i\mid X_1,\ldots, X_{i-1}).\]
\end{proposition}

The subadditivity can be extended to the setting of conditional entropy too.
Together with the chain rule, it shows that entropy can only go up when a condition is dropped.

\begin{proposition}
    For any three random variables $X,Y,Z$ with finite supports,
    \[H(X,Y\mid Z)\leq H(X\mid Z)+H(Y\mid Z).\]
    Equivalently, $H(X\mid Y,Z)\leq H(X\mid Z)$.
\end{proposition}

Lastly, we will need the following identity about adding a random variable that is determined by the others.
\begin{proposition}\label{proposition:NoMoreInfo}
    For any three random variables $X,Y,Z$ with finite supports such that $Z$ is determined by $X,Y$, i.e. $Z=f(X,Y)$ for some deterministic function $f$, 
    \[H(X,Z\mid Y)=H(X\mid Y).\]
\end{proposition}

\section{Rainbow triangle}\label{section:RainbowTriangle}
In this section, we first give a full proof of \cref{theorem:RainbowTriangle}, which will be followed by a discussion on its natural generalization.

\subsection{Proof of \cref{theorem:RainbowTriangle}}
As a reminder, we state the theorem here again.
\begin{theorem11}
Let $G=(V,E)$ be a simple graph and let edges be colored in colors red, green, and blue, so that each edge receives exactly one color. Let $R,G,B$ denote the number of red, green, blue edges, respectively. Let $T$ be the number of rainbow triangles in the graph.
Then $T^2\leq 2RGB$. 
\end{theorem11}

Before diving into the proof, we would like to highlight the importance of the assumption that $G$ is simple. 
When $G$ is allowed to be a multigraph, we may consider the following construction: in a graph on $N$ vertices, put an edge between any two vertices for each color.
Then we get $(1-o(1))N^3$ rainbow triangles with $(1-o(1))N^2/2$ edges of each color, showing that $T^2/RGB$ can be arbitrarily close to $8$.
This urges us to use the condition that $G$ is simple in a crucial way.

The idea of the proof is to understand the loss in a na\"ive application of Shearer's inequality. 
Suppose we pick a rainbow triangle chosen uniformly at random with vertices $X_1,X_2,X_3$ where $X_1X_2$ is red, $X_2X_3$ is green, and $X_1X_3$ is blue.
Shearer's inequality gives
\[2\log_2 T=2H(X_1,X_2,X_3)\leq H(X_1,X_2)+H(X_2,X_3)+H(X_1,X_3)\leq \log_2 (2R)+\log_2 (2G)+\log_2 (2B),\]
i.e. $T^2\leq 8RGB$. The factors of two come from the ``orientations'' of the edges, i.e. $(u,v)$ and $(v,u)$ could both be in the support of $(X_1,X_2)$. 
However, in the tight example (the $K_4$ blown-up, \cref{fig:rt}), Shearer's inequality is not tight as there are fewer choices for $X_2$ when $X_1$ is specified.
One might hope that we could fix this by improving the step where we apply Shearer's inequality, but that step could be tight when $G$ is a complete tripartite graph that is a blowup of a rainbow triangle.
In this particular example, the inequality that is not tight is the last one, as given a red edge $\{X_1,X_2\}$ we can uniquely determine $X_1$ and $X_2$.
In addition, although in both examples we have $H(\{X_1,X_2\})=\log_2R$, it is not necessarily the case when the graph is not as symmetric.
Thus, our strategy is to re-sample some of the random variables to avoid the aforementioned loss and treat carefully the conditional entropies that would appear in the proof of Shearer's inequality. 
\begin{proof}
In this proof, we consider the following random variables. Let $\Delta$ be a random rainbow triangle chosen uniformly at random, and let $\ell_r,\ell_g,\ell_b$ be its red, green, blue edges respectively. We emphasize here that $\ell_r,\ell_g,\ell_b$ and the edges we will consider in the proof are unordered pairs, instead of being ordered pairs as in the discussion above. Also denote by $v_r$ the vertex in $\Delta$ opposite to the edge $\ell_r$, and define $v_g$ and $v_b$ similarly.
For any vertex $v$, let $N_g(v)$ and $N_b(v)$ be the sets of the green neighbors and the blue neighbors of $v$ respectively. Let $u_b$ and $u_g$ be two random vertices sampled uniformly from $N_g(v_r)$ and $N_b(v_r)$ respectively, such that $u_{b},u_{g},\ell_r$ are conditionally independent given $v_{r}$. Let $L_g=\{v_{r},u_{b}\}, L_b=\{v_{r},u_{g}\}$ be unordered pairs as well.
Finally, sample $L^*_r$ from all the red edges uniformly at random and independently from all the other random variables. The random variables sampled above are shown in \cref{fig:rv}.

\begin{figure}[h]
    \centering
    
    \begin{tikzpicture}
    \coordinate (b) at (0,0);
    \coordinate (g) at (2,0);
    \coordinate (r) at (1,{sqrt(3)});
    \coordinate (B) at (-0.648776,0.6);
    \coordinate (G) at (2.648776,0.6);
    
    \coordinate (lr) at (1,0);
    \coordinate (lg) at (0.5,{sqrt(3)/2});
    \coordinate (lb) at (1.5,{sqrt(3)/2});
    
    \coordinate (lB) at ({(-0.648776+1)/2},{(0.6+sqrt(3))/2});
    \coordinate (lG) at ({(2.648776+1)/2},{(0.6+sqrt(3))/2});

    \coordinate (x) at (0,-1);
    \coordinate (y) at (2,-1);
    \coordinate (lR) at (1,-1);
    \draw[green, line width=1.4pt] (r) -- (b);
    \draw[red, line width=1.4pt] (b) -- (g);
    \draw[blue, line width=1.4pt] (g) -- (r);
    
    \draw[blue, line width=1.4pt] (r) -- (G);
    \draw[green, line width=1.4pt] (r) -- (B);

    \draw[red, line width=1.4pt] (x) -- (y);
    \node[above] at (r) {$v_r$};
    \node[below] at (g) {$v_g$};
    \node[below] at (b) {$v_b$};
    
    \node[above] at (lr) {$\ell_r$};
    \node[below,xshift=7,yshift=7] at (lg) {$\ell_g$};
    \node[below,xshift=-7,yshift=7] at (lb) {$\ell_b$};
    
    \node[above,yshift=1] at (B) {$u_b$};
    \node[above,yshift=1] at (G) {$u_g$};
    
    \node[above,xshift=2.5,yshift=1] at (lG) {$L_b$};
    \node[above,xshift=-2.5,yshift=-1] at (lB) {$L_g$};
    
    \node[below] at (lR) {$L^*_r$};

    \draw [fill] (r) circle (1.6pt);
    \draw [fill] (g) circle (1.6pt);
    \draw [fill] (b) circle (1.6pt);
    \draw [fill] (G) circle (1.6pt);
    \draw [fill] (B) circle (1.6pt);
    \draw [fill] (x) circle (1.6pt);
    \draw [fill] (y) circle (1.6pt);
\end{tikzpicture}
    \caption{The random variables sampled in the proof of \cref{theorem:RainbowTriangle}.}
    \label{fig:rv}
\end{figure}
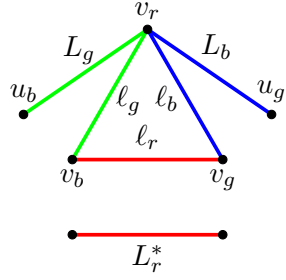

For convenience, we write $v\sim \ell$ for some $v\in V, \ell\in E$ if $v$ is the opposite vertex of $\ell$ in a rainbow triangle. 
In this case, we say that $v$ and $\ell$ \emph{form a rainbow triangle}.
Also, for any $u,v\in V$, we write $u\sim_r v$ if $\{u,v\}$ is a red edge.

Now our goal is to compute the entropy $2H(\Delta)-H(L^*_r)-H(L_g)-H(L_b)$. 
First, note that, by the chain rule, $H(\Delta)=H(\ell_r\mid v_{r})+H(v_{r})=H(\ell_r\mid v_{r})-H(L^*_r\mid v_{r})+H(L^*_r)+H(v_{r})$, where we are using the fact $H(L^*_r\mid v_r)=H(L^*_r)$ since $L^*_r$ and $v_{r}$ are independent. 
Conditioning on $v_{r}$, the entropy of $\ell_r$ is the logarithm of the number of red edges that forms a triangle with $v_{r}$, and the entropy of $L^*_r$ is the logarithm of the number of red edges. 
Therefore, we have $H(\ell_r\mid v_r)-H(L^*_r\mid v_r)=\mathbb{E}\left[\log_2 \mathbb{P}(v_{r}\sim L^*_r\mid v_{r})\right]$, and hence
\begin{align}\label{entropyR}
    H(\Delta)=H(L^*_r)+H(v_{r})+\mathbb{E}\left[\log_2 \mathbb{P}(v_{r}\sim L^*_r\mid v_{r})\right].
\end{align}
Next, by the chain rule, we have $H(\Delta)=H(v_{r})+H(v_{b},v_{g}\mid v_{r})$, and
\[H(v_{r})=H(v_r,L_g)-H(L_g\mid v_{r})=H(L_g)+H(v_r\mid L_g)-H(u_b\mid v_r).\]
Thus,
\begin{align}\label{entropyG}
    H(\Delta)=H(L_g)+H(v_{r}\mid L_g)-H(u_{b}\mid v_{r})+H(v_{b},v_{g}\mid v_{r}).
\end{align}
Similarly, we also have
\begin{align}\label{entropyB}
    H(\Delta)=H(L_b)+H(v_{r}\mid L_b)-H(u_{g}\mid v_{r})+H(v_{b},v_{g}\mid v_{r}).
\end{align}
We will sum up \cref{entropyR,entropyG,entropyB} and certain terms will cancel out.
To see this, note that
\begin{align}\label{entropy1}
H(v_{b},v_{g}\mid v_{r})+H(v_{r})=H(\Delta),
\end{align}
and also, similar as before, conditioning on $v_{r}$, the entropy $H(v_{b},v_{g}\mid v_{r})$ is equal to the logarithm of the number of red edges that form rainbow triangles with $v_r$. As $u_b$ and $u_g$ are uniformly chosen in $N_g(v_r)$ and $N_b(v_r)$, respectively, we have
\begin{align}\label{entropy2}
H(v_{b},v_{g}\mid v_{r})-H(u_{b}\mid v_{r})-H(u_{g}\mid v_{r})=\mathbb{E}\left[\log_2\mathbb{P}(u_{b}\sim_r u_{g}\mid v_{r})\right].
\end{align}
Thus, by adding \cref{entropyR,entropyG,entropyB,entropy1,entropy2}, we have
\begin{align}
    &2H(\Delta)-H(L^*_r)-H(L_g)-H(L_b)\nonumber\\
    =&\mathbb{E}\left[\log_2 \mathbb{P}(v_{r}\sim L^*_r\mid v_{r})\right]+\mathbb{E}\left[\log_2\mathbb{P}(u_{b}\sim_r u_{g}\mid v_{r})\right]+H(v_{r}\mid L_g)+H(v_{r}\mid L_b).\label{entropyALL}
\end{align}
Our goal now is to show that the right hand side is at most $1$. Note that, at this point, we already have the theorem with constant being $4$ instead of $2$ because the first two terms are both at most zero and the last two terms are both at most $1$, and hence $2\log_2 T\leq \log_2 R+\log_2 G+\log_2 B+2$.

For any $v\in V$, let $T_v$ denote the number of rainbow triangles such that $v$ is the vertex opposite to the red edge in this triangle, and let $\deg_g(v)\eqdef\abs{N_g(v)},\deg_b(v)\eqdef\abs{N_b(v)}$.
Then $\mathbb{P}(v_{r}=x)=T_x/T$. Thus, we have $\mathbb{P}(L_g=\{x,y\})=w_x+w_y$ where $w_x\eqdef\frac{T_x}{T\deg_g(x)}= \mathbb{P}(v_{r}=x,L_g=\{x,y\})$ and similar for $w_y$. 
As a consequence, for any green edge $\{x,y\}$,
\[H(v_{r}\mid L_g=\{x,y\})=-\frac{w_x}{w_x+w_y}\log_2\left(\frac{w_x}{w_x+w_y}\right)-\frac{w_y}{w_x+w_y}\log_2\left(\frac{w_y}{w_x+w_y}\right).\]
With foresight, we would like to show that
\[H(v_{r}\mid L_g=\{x,y\})\leq -\frac{w_x}{w_x+w_y}\log_2\left(\frac{T_x}{T_x+T_y}\right)-\frac{w_y}{w_x+w_y}\log_2\left(\frac{T_y}{T_x+T_y}\right),\]
which follows from the following more general claim.
\begin{claim}
For any fixed $c\in [0,1]$, the minimum value of $-c\log_2 z-(1-c)\log_2(1-z)$ for $z\in [0,1]$ is $h(c)\eqdef-c\log_2 c-(1-c)\log_2 (1-c)$.
\end{claim}
\begin{proof}
    The derivative of $-c\log_2 z-(1-c)\log_2(1-z)$ is $(-c/z+(1-c)/(1-z))\cdot \log_2 e$, which is positive if $z>c$ and negative if $z<c$.
    Thus, the minimum occurs at $z=c$.
\end{proof}
Therefore, by taking the expectation over all the green edges $\{x,y\}$, we have
\begin{align}
    H(v_{r}\mid L_g)&\leq \sum_{\{x,y\}\textup{ is green}}\mathbb{P}\left(L_g=\{x,y\}\right)\left(-\frac{w_x}{w_x+w_y}\log_2\left(\frac{T_x}{T_x+T_y}\right)-\frac{w_y}{w_x+w_y}\log_2\left(\frac{T_y}{T_x+T_y}\right)\right)\nonumber\\
    &=\sum_{\{x,y\}\textup{ is green}}-w_x\log_2\left(\frac{T_x}{T_x+T_y}\right)-w_y\log_2\left(\frac{T_y}{T_x+T_y}\right)\nonumber\\
    &=\sum_x\sum_{y: \{x,y\}\textup{ is green}}-w_x\log_2\left(\frac{T_x}{T_x+T_y}\right)\nonumber\\
    &=\sum_{x}\frac{T_x}{T}\sum_{y:\{x,y\}\textup{ is green}}-\frac{1}{\deg_g (x)}\log_2\left(\frac{T_x}{T_x+T_y}\right)\nonumber\\
    &=\mathbb{E}\left[\mathbb{E}\left[-\log_2\left(\frac{T_{v_{r}}}{T_{v_{r}}+T_{u_{b}}}\right)\,\middle|\, v_{r}\right]\right]\nonumber\\
    \label{entropyAfterAll}&\leq \mathbb{E}\left[-\log_2\left(\frac{T_{v_{r}}}{T_{v_{r}}+\mathbb{E}[T_{u_{b}}\mid v_{r}]}\right)\right],
\end{align}
where the last inequality follows from the concavity of logarithm.
We can bound this term by the first two terms in the right hand side of \cref{entropyALL} via the following claim.
\begin{claim}\label{claim:combinatorics}
    We have 
    \[\mathbb{P}(v_{r}\sim L^*_r\mid v_{r})\mathbb{P}(u_{b}\sim_r u_{g}\mid v_{r})\leq \frac{T_{v_{r}}}{T_{v_{r}}+\mathbb{E}[T_{u_{b}}\mid v_{r}]}.\]
\end{claim}
\begin{proof}
    For each $x\in N_g(v_r)$, we make the following notations.
    For each rainbow triangle counted in $T_{x}$, its red edge either forms a rainbow triangle with $v_{r}$ or not. We call the first kind of rainbow triangle type A and the second type B, and also denote the number of type A and type B triangles by $A_{v_{r},x}$ and $B_{v_{r},x}$ respectively.
    By definition, $T_{x}=A_{v_{r},x}+B_{v_{r},x}$.

    For each type A triangle, there is either a blue edge or a green edge connecting $N_g(v_{r})$ and $N_b(v_{r})$. 
    For each $x\in N_g(v_r)$, summing over all type $A$ triangles counted in $T_{x}$, we count $A_{v_r,x}$ non-red edges between $x$ and $N_b(v_r)$ with multiplicities.
    Therefore when we sum over all $x\in N_g(v_r)$, we have counted $\sum_{x\in N_{g}(v_r)}A_{v_{r},x}$ non-red edges between $N_g(v_r)$ and $N_b(v_r)$ with multiplicities, where the multiplicities are at most $\deg_g(v_{r})-1$.
    Thus, there are in total at least 
    \[\frac{\sum_{x\in N_{g}(v_r)}A_{v_{r},x}}{\deg_g(v_{r})-1}> \mathbb{E}[A_{v_{r},u_{b}}\mid v_{r}]\]
    such edges. As a consequence,
    \[\mathbb{P}(u_{b}\sim_r u_{g}\mid v_{r})\leq \frac{T_{v_{r}}}{T_{v_{r}}+\mathbb{E}[A_{v_{r},u_{b}}\mid v_{r}]}.\]

    For each type B triangle, it gives a red edge that does not form a rainbow triangle with $v_{r}$, and there are at least $\max_{u_{b}}B_{v_{r},u_{b}}$ many such edges. Thus,
    \[\mathbb{P}(v_{r}\sim L^*_r\mid v_{r})\leq \frac{T_{v_{r}}}{T_{v_{r}}+\max_{u_{b}}B_{v_{r},u_{b}}}\leq \frac{T_{v_{r}}}{T_{v_{r}}+\mathbb{E}[B_{v_{r},u_{b}}\mid v_{r}]}.\]

    Combining these two inequalities, we have
    \[\mathbb{P}(v_{r}\sim L^*_r\mid v_{r})\mathbb{P}(u_{b}\sim_r u_{g}\mid v_{r})\leq \frac{T_{v_{r}}}{T_{v_{r}}+\mathbb{E}[B_{v_{r},u_{b}}\mid v_{r}]}\frac{T_{v_{r}}}{T_{v_{r}}+\mathbb{E}[A_{v_{r},u_{b}}\mid v_{r}]}\leq \frac{T_{v_{r}}}{T_{v_{r}}+\mathbb{E}[T_{u_{b}}\mid v_{r}]}.\qedhere\]
\end{proof}
Combining \cref{entropyALL}, \cref{entropyAfterAll} and the claim above, we get
\[2H(\Delta)-H(L^*_r)-H(L_g)-H(L_b)\leq H(v_{r}\mid L_b)\leq 1.\]
Since $H(\Delta)=\log_2 T$, $H(L^*_r)\leq \log_2 R$, $H(L_g)\leq \log_2 G$, and $H(L_b)\leq \log_2 B$ by \cref{proposition:UnifBound}, we have $T^2\leq 2RGB$.
\end{proof}

\subsection{Higher-uniformity analogue}
Inspired by the joints problem, the natural generalization of the rainbow triangle problem is the following.
\begin{question}
    Let $d\geq 2$ be a positive integer, and let $c_1,\ldots, c_d$ be $d$ colors.
    In a simple $(d-1)$-uniform hypergraph $H$ colored by the $d$ colors, let $C_i$ be the number of hyperedges with the color $c_i$ for any $i\in[d]$.
    Here, a hypergraph is simple if no hyperedges appear twice.
    Let $T$ be the number of rainbow $K^{(d-1)}_d$'s in $H$.
    How large can the ratio $T^{d-1}/C_1\cdots C_d$ be?
\end{question}

We will denote by $\kappa_d$ the supremum of the possible values for the ratio given a fixed $d$.
Note that it is always finite by Shearer's inequality as follows.
We first sample uniformly a rainbow $d$-clique $\Delta=\{X_1,\ldots, X_d\}$ from $H$, where $\Delta\backslash \{X_i\}$ has the $i$-th color.
Then
\begin{align*}
  (d-1)H(\Delta) &= (d-1)H(X_1,\ldots, X_d)\\
  &\leq \sum_{i=1}^{d}H(X_1,\ldots, \widehat{X_i},\ldots, X_d)\\
&\leq \sum_{i=1}^{d}\left(\log_2 (d-1)!+H(\Delta\backslash \{X_i\})\right) \\
&\leq d\log_2(d-1)!+\sum_{i=1}^{d}\log_2C_i.
\end{align*}
As $H(\Delta)=\log_2T$, by rearranging, we immediately get $\kappa_d\leq \left((d-1)!\right)^d$.

In fact, part of the proof of \cref{theorem:RainbowTriangle} generalizes to give that $\kappa_d\leq \frac{1}{2}\prod_{i=1}^{d-1}i^i$ whenever $d\geq 3$.
This can be done by an induction on $d$.
The $d=3$ case was already done, so now let us assume that it holds for $d-1$.
We again sample a rainbow $d$-clique $\Delta=\{X_1,\ldots, X_d\}$ in $H$ uniformly at random.
For any $i\in[d-1]$, re-sample a hyperedge $L_i$ with color $c_i$ that contains $X_d$ uniformly, where it is done conditionally indenpendent of $X_1,\ldots, X_{d-1}$.
Then for each $i\in[d-1]$,
\[H(\Delta) = H(X_d)+H(\Delta\backslash \{X_d\}\mid X_d) = H(L_i)+H(X_d\mid L_i)-H(L_i\mid X_d)+H(\Delta\backslash \{X_d\}\mid X_d).\]
Note that $H(X_d\mid L_i)\leq \log_2(d-1)$ and $H(L_i)\leq \log_2C_i$ for all $i\in[d-1]$.
Adding all $d-1$ equations up, moving terms around, and upper bounding one of the terms $H(\Delta\backslash \{X_d\}\mid X_d)$ by $\log_2C_d$, we get
\[(d-1)H(\Delta)-\sum_{i=1}^{d}\log_2C_i\leq (d-1)\log_2(d-1)+(d-2)H(\Delta\backslash \{X_d\}\mid X_d)-\sum_{i=1}^{d-1}H(L_i\mid X_d).\]
Now for every $v\in V(H)$, let $H/v$ be the $(d-2)$-uniform hypergraph on $V(H)\backslash \{v\}$ where $e\in E(H/v)$ if and only if $e\cup\{v\}\in E(H)$, and they are colored the same if it is the case.
Note that $H(L_i\mid X_d=v)$ is exactly the logarithm of the number of $(d-2)$-edges in $H/v$ with the color $c_i$ as $L_i$ is re-sampled uniformly.
In addition, when $X_d=v$, $\Delta\backslash \{X_d\}$ is a rainbow $(d-1)$-clique avoiding the $d$-th color in $H/v$, and thus $H(\Delta\backslash \{X_d\}\mid X_d=v)\leq \log_2T'_v$ where $T'_v$ is the number of rainbow $(d-1)$-cliques in $H/v$ avoiding the $d$-th color.
Thus, by the inductive hypothesis,
\[(d-2)H(\Delta\backslash \{X_d\}\mid X_d)-\sum_{i=1}^{d-1}H(L_i\mid X_d)\leq \log_2\left(\frac{1}{2}\prod_{i=1}^{d-2}i^i\right),\]
and hence the inductive step is done.

Although this bound is tight for $d=3$, it is never tight when $d>3$.
Indeed, by directly applying the joints bound, we get the following stronger upper bound of $\kappa_d$.
\begin{theorem}
    For any $d\geq 2$, we have $\kappa_d\leq d!$.
\end{theorem}
This is an immediate corollary of Theorem 3.2 in \cite{YZ23} once applied to the generically induced configurations.

One of the reasons why the inductive steps are not tight is that we have no higher uniformity analog of \cref{claim:combinatorics}.
We also could not find examples where $H(X_d\mid L_i)\leq \log_2(d-1)$ is tight for lots of $i$'s, though we currently have no ways of bounding this more cleverly.

For the lower bound, we can generalize the construction for $d=3$ to odd $d$.
For even $d$'s, we will use the following monotonicity of $\kappa_d$.

\begin{proposition}\label{proposition:MonotoneKappa}
    For $d\geq 3$, we have $\kappa_{d+1}\geq \kappa_d$.
\end{proposition}
\begin{proof}
Let $H$ be a $(d-1)$-uniform hypergraph colored by $d$ colors. 
Assume there are $C_i$ hyperedges with color $c_i$ for $i\in [d]$ and there are $T$ rainbow $K_d^{(d-1)}$'s.
We construct a $(d+1)$-colored $d$-uniform hypergraph $H'$ by adding a new distinguished vertex $v$ into all the edges, and color all the rainbow $(d-1)$-clique in the $d$-th color.
This way, any rainbow $(d-1)$-clique is extended to a rainbow $d$-clique by adding $v$ in. It follows that the number $C'_i$ of $c_i$-colored hyperedges equals $C_i$ for $i\in [d]$, the number $C'_{d+1}$ of $c_{d+1}$-colored hyperedges equals $T$, and the number $T'$ of rainbow $K_{d+1}^{(d)}$ in $H'$ also equals $T$.
Therefore, the ratio $T'^{d}/C'_1\dots C'_{d+1}$ equals the ratio $T^{d-1}/C_1\dots C_{d}$, showing that $\kappa_{d+1}\geq \kappa_d$.
\end{proof}

\begin{theorem}
    For any $d\geq 2$, we have
    \[
    \kappa_d\geq \begin{cases}
        2^d/(d+1)&\textup{ if }d\textup{ is odd;}\\
        2^{d-1}/d&\textup{ if }d\textup{ is even.}
    \end{cases}
    \]
\end{theorem}
\begin{proof}
    For odd $d$, consider any decomposition of $K_{d+1}$ into $d$ matchings $M_1,\ldots, M_d$.
    Now for each hyperedge in $K^{(d-1)}_{d+1}$, color it with the $i$-th color if its complement is in $M_i$.
    Then any $d$-clique in $K^{(d-1)}_{d+1}$ is rainbow, so we get $T=d+1$ and $C_1=\cdots=C_d=(d+1)/2$.
    Thus,
    \[\kappa_d\geq \frac{(d+1)^{d-1}}{\left(\frac{d+1}{2}\right)^d}=\frac{2^d}{d+1}.\]

    For $d=2$, it is clear that $\kappa_d=1$.
    For other even $d$'s, by \cref{proposition:MonotoneKappa} we have 
    \[\kappa_d\geq\kappa_{d-1}\geq \frac{2^{d-1}}{d}.\qedhere\]
\end{proof}

For $d=4$, we have a better construction than the one above, where the constant we get is $3$.
As the construction is harder to state and there is no reason to believe that it is optimal, we include the construction in \cref{appendix:tetrahedra}.

\section{A generalization of Kruskal--Katona}\label{section:KruskalKatona}
In this section, we first give an entropy-based proof for the Lov\'asz's version of the Kruskal--Katona theorem.
Then we give a generalization of the proof that has more applications.
\subsection{Entropy-based Proof for Kruskal--Katona}\label{subsection:SimpleKK}

\begin{theorem12}
Let $\cA$ be a family of $d$-subsets of $[n]$ and $\partial \cA$ be all the $(d-1)$-subsets that are contained in some sets in $\cA$. If $\abs{\cA}=\binom{t}{d}$ for some real number $t\geq d$, then $\abs{\partial\cA}\geq \binom{t}{d-1}$.
\end{theorem12}
\begin{proof}
    Pick a random set $A\in\cA$ uniformly at random, and then pick a random ordering $X_1,\dots,X_d$ of elements in $A$ uniformly at random. It follows that \[H(X_1,\dots, X_d)=\log_2 (d!|\cA|)=\sum_{k=1}^d \log_2(t+1-k).\]
    We also have
\begin{align*}
H(X_1,\ldots,X_d)&=H(X_1,\ldots,X_{d-1})+H(X_d\mid X_1,\ldots,X_{d-1}).
\end{align*}
Assume $H(X_d\mid X_1,\ldots,X_{d-1})=\log_2 r$. 
We will prove that for any $1\leq k\leq d-1$, 
\begin{align*}
    H(X_k\mid X_1,\ldots,X_{k-1})\geq \log_2 (r+d-k).
\end{align*}
Once we prove this inequality, we have
\[H(X_1,\ldots,X_d)=\sum_{k=1}^d H(X_k\mid X_1,\ldots,X_{k-1})\geq \sum_{k=1}^d\log_2 (r+d-k),\]
so $r\leq t-d+1$. Since $\{X_1,\ldots,X_{d-1}\}$ is a set in $\partial\cA$, we have
\begin{align*}
\log_2 t(t-1)\ldots (t-d+1)&=H(X_1,\ldots,X_d)\\
&=H(X_1,\ldots,X_{d-1})+H(X_d\mid X_1,\ldots,X_{d-1})\\
&=H(X_1,\ldots,X_{d-1})+\log_2 r\\
&\leq \log_2 \left((d-1)!|\partial\cA|\right)+\log_2 (t-d+1).
\end{align*}
By rearranging the terms, we get
\[|\partial\cA|\geq \frac{t(t-1)\ldots(t-d+2)}{(d-1)!}=\binom{t}{d-1}.\]
It remains to prove the inequality. Indeed, we will prove that for each $k=1,\ldots, d-1$,
\begin{align}\label{equation:key}
    2^{H(X_k\mid X_1,\ldots,X_{k-1})}\geq 2^{H(X_{k+1}\mid X_1,\ldots,X_{k})}+1,
\end{align}
which implies the desired inequality via a telescoping sum.
We will fix a $k$ in this proof from now on.

First flip a coin $\coin \sim \text{Bern}(p)$ that is independent from $X_1,\dots,X_d$, where $p\in [0,1]$ is a number to be determined later. Set $X^*=X_{k+1}$ if $\coin =1$. Otherwise, set $X^*=X_{k}$. Note that, when conditioning on $X_1,\dots,X_{k-1}$, the distributions of $X_k,X_{k+1}$ and $X^*$ are identical. Therefore,
\begin{align*}
H(X_k\mid X_1,\dots,X_{k-1})=&H(X^*\mid X_1,\dots,X_{k-1})\\
\geq &H(X^*\mid X_1,\dots,X_{k-1},X_k)
\end{align*}
Moreover, conditioning on $X_1=x_1,\dots,X_{k}=x_k$ for any $x_1,\dots,x_k$ in the support of $X_1,\dots,X_k$, the support of $X_{k+1}$ does not contain $x_k$, and hence $\coin$ is determined by $X_1,\dots,X_k,X^*$.
Thus, by \cref{proposition:NoMoreInfo} we have
\begin{align*}
 &H(X^*\mid X_1,\dots,X_{k})\\
=&H(X^*,\coin \mid X_1,\dots,X_{k})\\
=&H(X^*\mid \coin ,X_1,\dots,X_{k})+H(\coin\mid X_1,\dots,X_{k})\\
=&p H(X_{k+1}\mid \coin=1,X_1,\dots,X_{k})+(1-p)H(X_{k}\mid \coin=0,X_1,\dots,X_{k})+h(p).\\
=&p H(X_{k+1}\mid X_1,\dots,X_{k})+h(p).
\end{align*}
By combining the inequalities above, we have
\[H(X_k|X_1,\dots,X_{k-1})\geq pH(X_{k+1}\mid X_1,\dots,X_k)+h(p).\]
Let $s=2^{H(X_{k+1}\mid X_1,\dots,X_{k})}$ and pick $p=\frac{s}{s+1}$. It follows that 
\begin{align*}
 H(X_k\mid X_1,\dots,X_{k-1})\geq & \frac{s}{s+1}\log_2 s-\frac{s}{s+1}\log_2 \left(\frac{s}{s+1}\right)-\frac{1}{s+1}\log_2 \left(\frac{1}{s+1}\right)\\
 =&\log_2(s+1).   
\end{align*}
Equivalently, we have 
\begin{align*}
    2^{H(X_k\mid X_1,\ldots,X_{k-1})}\geq 2^{H(X_{k+1}\mid X_1,\ldots,X_{k})}+1.
\end{align*}
This completes the proof.
\end{proof}
\subsection{Generalized Kruskal--Katona}
In this section, we first state the necessary definitions in order to state the main result \cref{theorem:entropykk}.
We then prove the main theorem and show \cref{theorem:Lovasz} and \cref{theorem:q-analog} as immediate corollaries.

\cref{definition:entropykk} describes the basic settings of our main theorem of this section. Instead of working with set systems, we will focus on tuples. 
In our setting, we consider a family $\cF$ consisting of some $d$-tuples and we are interested in its ``$(d-1)$-shadow'', which is the set of all the $(d-1)$-tuples obtained from dropping the last entry.
The family $\cF$ cannot be arbitrary, which will be made precise by the definition and the theorem statement.

In the proof in \cref{subsection:SimpleKK}, we saw that $2^{H(X_{k+1}\mid X_1,\ldots,X_{k-1})}\geq 2^{H(X_{k+1}\mid X_1,\ldots,X_{k})}+1$. 
We can think of this inequality as saying that, once we reveal $X_k$, we have one less possible outcome for $X_{k+1}$ to pick from, so we can get one more than the na\"ive bound given by dropping the conditioning random variable $2^{H(X_{k+1}\mid X_1,\ldots,X_{k-1})}\geq 2^{H(X_{k+1}\mid X_1,\ldots,X_{k})}$. Moreover, we can get $2^{H(X_k)}\geq 2^{H(X_k\mid X_1,\ldots,X_{k-1})}+k-1$ from the telescoping sum. This inequality can be interpreted as follows: revealing $X_1,\dots,X_{k-1}$ forbids $c_{k-1}=k-1$ possible outcomes of $X_k$, so we gain a $k-1$ on the right hand side. We can generalize this idea and define \emph{$(c_1,\dots,c_{d-1})$-forbidding systems}.

In the first part of the definition, we define certain conditions we need for good and bad families of tuples.
Basically, for any fixed $x_1,\ldots,x_{k-1}$ with $(x_1,\ldots,x_{k-1},*,\ldots,*)$ being good, the number of $x_k$'s such that $(x_1,\ldots,x_{k-1},x_{k},*,\ldots,*)$ is bad is exactly $c_k$.
In the actual definition, we will consider good and bad multisets instead in order to always keep the families symmetric, and we think of a tuple as bad if some of its entries form a bad multiset.
In the second part, we define the notion of compatible set of elements.
Roughly speaking, a compatible set is where the structure of the good and bad families restrict nicely to.
Finally, our family $\mathcal{F}$ must consist of all the good $d$-tuples from several compatible sets. For example, in the case of the Kruskal--Katona theorem, a compatible set is just a set $A$ of $d$ distinct elements, and we picked the $d$-tuples from all the possible ordering of elements in $A\in\cA$.

Throughout this section, we will denote by $B$ the collection of bad objects and denote by $G$ the collection of good objects.

\begin{definition}\label{definition:entropykk}
Given a set $U$, let $M_k$ be the family of all the multisets of size $k$ with elements from $U$. In the following, we assume all the set operators are for multisets.

Let $d>0$ be an integer. Suppose $\mathbf{B}=(\emptyset=B_1,B_2,\ldots,B_d)$ and $\mathbf{G}=(G_1,G_2,\ldots,G_d)$ are two sequences of families of multisets with $B_k\sqcup G_k= M_k$ for all $k$. For any $k$, any multiset $A\in M_k$, and any $x\in U$, we call $A\cup\{x\}$ an \emph{extension} of $A$, and say that $A\cup \{x\}$ is a bad (good) extension if $A\cup\{x\}\in B_{k+1}$ ($G_{k+1}$ resp.). Let $0\leq c_1\leq c_2\leq \ldots\leq c_{d-1}$ be integers. We say that $(\mathbf{B},\mathbf{G})$ is a $(c_1,c_2,\ldots,c_{d-1})$-\emph{forbidding system} if the following holds. For any multisets $A\in M_k$ with $0<k<d$, if $A\in B_k$, then all the extensions of $A$ are bad. Otherwise, there are exactly $c_{k}$ bad extensions of $A$. 

Let $S$ be a subset of $U$.
We say that $S$ is \emph{$(\mathbf{B},\mathbf{G})$-compatible} if the following holds. For all $0<k<d$ and 
all multisets $A\in G_k$ that consist of elements from $S$, all $x\in U$ such that $A\cup\{x\}\in B_{k+1}$ are in $S$. If $S$ is $(\mathbf{B},\mathbf{G})$-compatible, then we write \[S^{(d)}\eqdef\{(x_1,\ldots,x_d)\in S^d\mid \{x_1,\ldots,x_d\}\in G_d\}.\]
\end{definition}

We give an example to illustrate the idea behind the definition. The following example corresponds to the original version of the Kruskal--Katona theorem.

\begin{example}\label{example:kk}
Let $U=[n]$. Denote by $B_k$ the family of multisets of size $k$ that contain repeated elements, and hence $G_k$ is the family of sets (without repeated elements) of size $k$. For any multiset $A$ with repeated elements, we know that any extension of $A$ also contains repeated elements. If $A\in G_k$, then there are exactly $k$ bad extensions of $A$. Namely, those extensions that extend by the elements in $A$ are bad. Thus, the pair $(\mathbf{B},\mathbf{G})$ is a $(1,2,\ldots,d-1)$-forbidding system.

The $d$-tuples in $S^{(d)}$ from $(\mathbf{B},\mathbf{G})$-compatible sets $S$ will be the choices we can put in our family $\cF$. In this case, $S$ is $(\mathbf{B},\mathbf{G})$-compatible if and only if $S$ is a set. In the setting of the Kruskal--Katona theorem, we only take $S$ to be sets of size $d$.
\end{example}

Here are some properties of $S^{(d)}$ that will be useful later in the main proof.

\begin{proposition}\label{proposition:Balanced}
    Given a set $U$ and a $(c_1,c_2,\ldots,c_{d-1})$-forbidding system $(\mathbf{B},\mathbf{G})$ on $U$. Suppose $S\subseteq U$ is a $(\mathbf{B},\mathbf{G})$-competible set. It follows that
    \begin{itemize}
        \item the set $S^{(d)}$ is symmetric, i.e. $(x_{\sigma(1)},\ldots,x_{\sigma(d)})\in S^{(d)}$ for all $(x_1,\ldots,x_d)\in S^{(d)}$ and all permutations $\sigma$ of $1,2,\ldots, d$;
        \item the size of $S^{(d)}$ is $|S|(|S|-c_1)\ldots(|S|-c_{d-1})$;
        \item for any $1\leq k\leq d$ and any fixed $\{x_1,\ldots,x_k\}\in G_k$ with $x_1,\ldots,x_k\in S$,  the number of tuples of the form $(x_1,\ldots,x_k,*,\ldots,*)$ in $S^{(d)}$ is $(|S|-c_{k})\ldots(|S|-c_{d-1})$.
    \end{itemize}
\end{proposition}
\begin{proof}
    The first part follows from the definition. The second part and the third part follow from the following fact: If $\{x_1,\ldots,x_k\}\in G_k$ and $x_1,\ldots,x_k\in S$, then there are $|S|-c_{k}$ good extensions of $\{x_1,\ldots,x_k\}$. This is because $S$ contains all the $x$'s such that $\{x_1,\ldots,x_k,x\}\in B_{k+1}$.
\end{proof}

Now, we are ready to state our generalization of the Kruskal--Katona theorem.

\begin{theorem}\label{theorem:entropykk}
Given a set $U$. Let $(\mathbf{B},\mathbf{G})$ be a $(c_1,c_2,\ldots,c_{d-1})$-forbidding system. Assume $S_1,\ldots,S_m$ are $(\mathbf{B},\mathbf{G})$-compatible such that $S_1^{(d)},\ldots,S_m^{(d)}$ are mutually disjoint. Let $\cF=\cup_i S_i^{(d)}$ and let $\partial\cF$ be the family of all $(d-1)$-tuples obtained from $\cF$ by removing the last entry from any tuples in $\cF$. It follows that, if $t\geq c_{d-1}$ is a real such that $|\cF|=t(t-c_1)\ldots(t-c_{d-1})$, then
\[|\partial\cF|\geq t(t-c_1)\ldots(t-c_{d-2}).\]
\end{theorem}

We will begin with proving the generalization of \cref{equation:key}.
We start with a lemma, whose proof should be reminiscent of the proof of \cref{equation:key}.

\begin{lemma}\label{lemma:KK}
Given two sets $U,V$. Let $D_1\sqcup D_2=U\times V$ be a partition. Sample three random variables $X_1,X_2,Y$ such that $(X_i,Y)\in D_i$ for $i=1,2$. If the law of $X_1$ and $X_2$ are the same, then
\[2^{H(X_1)}\geq 2^{H(X_1\mid Y)}+2^{H(X_2\mid Y)}.\]
\end{lemma}
\begin{proof}
Let $p\in [0,1]$ be a number to be determined. Sample $X^*$ in the following way. First flip a coin $\coin \sim \text{Bern}(p)$ independently from $X_1,X_2,Y$. Set $X^*=X_1$ if $\coin =1$. Otherwise, set $X^*=X_2$. Note that $X^*$ also has the same law as $X_1$ and $X_2$. It follows that $H(X_1)=H(X^*)\geq H(X^*\mid Y)$. Observe that, given $Y=y$, the support of $X_1$ and $X_2$ are disjoint, and hence $\coin$ is determined by $X^*,Y$
Thus, by \cref{proposition:NoMoreInfo} we have
\begin{align*}
H(X^*\mid Y)&=H(X^*,\coin \mid Y)\\
&=H(X^*\mid \coin ,Y)+H(\coin\mid Y)\\
&=p H(X_1\mid Y)+(1-p)H(X_2\mid Y)+h(p).
\end{align*}
Therefore,
\begin{align*}
H(X_1)\geq p H(X_1\mid Y)+(1-p)H(X_2\mid Y)+h(p).
\end{align*}
Let $s_i=2^{H(X_i\mid Y)}$ for $i=1,2$ and pick $p=\frac{s_1}{s_1+s_2}$. It follows that \begin{align*}
    H(X_1)\geq& \frac{s_1}{s_1+s_2}\log_2 s_1+\frac{s_2}{s_1+s_2}\log_2 s_2-\frac{s_1}{s_1+s_2}\log_2 \left(\frac{s_1}{s_1+s_2}\right)-\frac{s_2}{s_1+s_2}\log_2 \left(\frac{s_2}{s_1+s_2}\right)\\
    =&\log_2 \left(s_1+s_2\right).\qedhere
\end{align*}
\end{proof}

We next extract the necessary definition from \cref{definition:entropykk}.
Recall that we would like to compare $2^{H(X_{k+1}\mid X_1,\ldots, X_{k-1})}$ with $2^{H(X_{k+1}\mid X_1,\ldots, X_{k})}$.
To do so, we will first fix $X_1,\ldots, X_{k-1}$ and focus on the joint distribution $(X_k, X_{k+1})$.
The properties we will need $(X_k,X_{k+1})$ to satisfy are the following.

\begin{definition}\label{definition:balanced}
Given two sets $U,V$ and a positive integer $c$. Let $B\sqcup G=U\times V$ be a partition. Again, we call a pair $(x,y)$ bad (good) if $(x,y)\in B$ ($(x,y)\in G$ resp.). Assume that, for any fixed $y\in V$, there are exactly $c$ pairs $(x,y)\in B$. In this case, we call the pair $(B,G)$ a \emph{$c$-regular} partition of $U\times V$.

We call a set $S\times T\subseteq U\times V$ \emph{balanced} with respect to $(B,G)$ if the following conditions hold. First, for any fixed $y\in T$, $S\times T$ contains all the bad $(x,y)$. Second, there exists an integer $r$ such that, for any fixed $x\in S$, the number of bad pairs $(x,y)$ in $S\times T$ is $r$.
\end{definition}

\begin{corollary}\label{cor:entropykk}
Given two sets $U,V$ and a positive integer $c$. Let $(B,G)$ be a $c$-regular partition of $U\times V$.
Let $\bS\times \bT$ be a random balanced set. Sample a random pair $(X,Y)$ from all the good pairs in $\bS\times \bT$ chosen uniformly at random. It follows that 
\[2^{H(X)}-2^{H(X\mid Y)}\geq c.\]
\end{corollary}
\begin{proof}
Let $X'$ be a random variable sampled uniformly at random from all the elements $x\in U$ with $(x,Y)$ being bad. Given any $Y$, we know that the are $c$ possible choices for $X'$. Therefore, $H(X'\mid Y)=\log_2 c$. Following from \cref{lemma:KK}, it suffices to check that $X$ and $X'$ have the same law.

Note that, for any balanced $S\times T$, the number of $y$ such that $(x,y)$ is good/bad is a constant for $x\in S$, and also the number of $x$ such that $(x,y)$ is good/bad is a constant for $y\in T$. 
Therefore, the probability $\mathbb{P}(X=x\mid \bS\times \bT=S\times T)$ and $\mathbb{P}(X'=x\mid \bS\times \bT=S\times T)$ are constants for $x\in S$, so they are both $1/|S|$. This implies that $X$ and $X'$ have the same law.
\end{proof}

Now, we are ready to proof the main theorem.

\begin{proof}[Proof of \cref{theorem:entropykk}]
Sample random variables $X_1,\ldots,X_d\in U$ in the following way. First, we sample a random index $i\in [m]$ uniformly at random, and then sample $(X_1,\ldots,X_d)\in S_i^{(d)}$ uniformly at random. Since $S_1^{(d)},\ldots,S_m^{(d)}$ are mutually disjoint, we have $H(X_1,\ldots,X_d)=\log_2 |\cF|$. We also have
\begin{align*}
H(X_1,\ldots,X_d)&=H(X_1,\ldots,X_{d-1})+H(X_d\mid X_1,\ldots,X_{d-1}).
\end{align*}
Assume $H(X_d\mid X_1,\ldots,X_{d-1})=\log_2 r$. 
We will prove that for any $1\leq k\leq d-1$, 
\[H(X_k\mid X_1,\ldots,X_{k-1})\geq \log_2 (r+c_{d-1}-c_{k-1})\] where $c_0=0$. 

We know that $\{X_1,\ldots,X_{k-1}\}\in G_{k-1}$. Fix any $\{x_1,\ldots,x_{k-1}\}\in G_{k-1}$. We will apply \cref{cor:entropykk} after conditioning on $X_1=x_1,\ldots,X_{k-1}=x_{k-1}$. Let $V=\{x\in U\mid \{x_1,\ldots,x_{k-1},x\}\in G_k\}$, which contains the set of possible outcomes for $X_k$. Define $B$ to be the set of bad pairs $(x,y)\in V^2$ with $\{x_1,\ldots,x_{k-1},x,y\}\in B_{k+1}$, and set $G=V^2\setminus B$. Note that, for any fixed $y\in V$, there are exactly $(c_k-c_{k-1})$ pairs $(x,y)\in B$. This is because there are $c_{k}$ choices of $x$ such that $\{x_1,\ldots,x_{k-1},x,y\}$ is bad and all but $c_{k-1}$ of which are in $V$. Thus, $(B,G)$ is a $(c_k-c_{k-1})$-regular partition of $V^2$.

For any $(\mathbf{B},\mathbf{G})$-compatible set $S$ that contains $x_1,\ldots,x_{k-1}$, we observe that $T\times T$ is balanced, where $T=S\cap V$. 
To see this, we first note that for any fixed $y\in T$, there are $(c_k-c_{k-1})$ ``bad pairs'' $(x,y)$ in $S\times S$ since all the possible $x$ such that $\{x_1,\ldots,x_{k-1},x,y\}\in B_{k+1}$ are in $S$.
Moreover, because of \cref{proposition:Balanced}, we know that
\[\mathbb{P}(X_1=x_1,\ldots,X_{k-1}=x_{k-1}\mid i=i_0)=\frac{1}{|S_{i_0}|(|S_{i_0}|-c_1)\dots(|S_{i_0}|-c_{k-2})}\]
holds for any $i_0$ such that $S_{i_0}$ contains $x_1,\ldots,x_{k-1}$, provided that $S_{i_0}^{(d)}$ is non-empty. Similarly, we know that
\[\mathbb{P}(X_1=x_1,\ldots,X_{k-1}=x_{k-1},X_k=x,X_{k+1}=y\mid i=i_0)=\frac{1}{|S_{i_0}|(|S_{i_0}|-c_1)\dots(|S_{i_0}|-c_{k})}\]
holds for any $i_0$ such that $S_{i_0}$ contains $x_1,\ldots,x_{k-1}$, provided that $S_{i_0}^{(d)}$ is non-empty, and any $(x,y)\in T_{i_0}^2\cap G$, where $T_{i_0}=S_{i_0}\cap V$. Thus, conditioning on $X_1=x_1,\ldots,X_{k-1}=x_{k-1}$, the random variables $X_{k},X_{k+1}$ are chosen in the same way as they were chosen from $T_i^2\cap G$ uniformly where $T_i=S_i\cap V$. 
Thus, we may apply \cref{cor:entropykk} and conclude that
\begin{align*}
    {H(X_{k+1}\mid X_1=x_1,\ldots,X_{k-1}=x_{k-1},X_k)}=&{H(X_{k}\mid X_1=x_1,\ldots,X_{k-1}=x_{k-1},X_{k+1})}\\
    \leq& \log_2\left(2^{H(X_k\mid X_1=x_1,\ldots,X_{k-1}=x_{k-1})}-(c_k-c_{k-1})\right).
\end{align*}
Note that $\log_2(2^s-c)$ is a concave function of $s$ for any $c\geq 0$. Thus, by taking the average over $(x_1,\dots,x_{k-1})$ (using the same law as $(X_1,\dots, X_{k-1})$) and using Jensen's inequality, we have
\[H(X_{k+1}\mid X_1,\ldots,X_{k-1},X_k)\leq \log_2\left(2^{H(X_k\mid X_1,\ldots,X_{k-1})}-(c_k-c_{k-1})\right)\]
holds for all $k=1,2,\ldots,d-1$. Thus, for each $k$, we have
\[2^{H(X_k\mid X_1,\ldots,X_{k-1})}\geq 2^{H(X_{k+1}\mid X_1,\ldots, X_k)}+(c_k-c_{k-1})\geq\cdots\geq 2^{H(X_d\mid X_1,\ldots,X_{d-1})}+(c_{d-1}-c_{k-1}),\]
showing that $H(X_k\mid X_1,\ldots,X_{k-1})\geq \log_2 (r+c_{d-1}-c_{k-1})$.

Summing it up, we have
\[H(X_1,\ldots,X_d)=\sum_{k=1}^d H(X_k\mid X_1,\ldots,X_{k-1})\geq \sum_{k=1}^d\log_2 (r+c_{d-1}-c_{k-1}),\]
so $t\geq r+c_{d-1}$. Since $(X_1,\ldots,X_{d-1})$ is supported on $\partial\cF$, we have
\begin{align*}
\log_2 t(t-c_1)\ldots (t-c_{d-1})&=H(X_1,\ldots,X_d)\\
&=H(X_1,\ldots,X_{d-1})+H(X_d\mid X_1,\ldots,X_{d-1})\\
&=H(X_1,\ldots,X_{d-1})+\log_2 r\\
&\leq \log_2 |\partial\cF|+\log_2 (t-c_{d-1}).
\end{align*}
By rearranging the terms, we get
\[|\partial\cF|\geq t(t-c_1)\ldots(t-c_{d-2}).\qedhere\]
\end{proof}

We now see how this general result implies Lov\'asz's version of the Kruskal--Katona theorem and its $q$-analog (\cref{theorem:Lovasz,theorem:q-analog}).

\begin{proof}[Proof of \cref{theorem:Lovasz} using \cref{theorem:entropykk}]
As we discussed in \cref{example:kk}, when $U=[n]$, $B_k$ is the family of multisets of size $k$ that contains repeated elements, and $G_k$ is the family of sets of size $k$, $(\mathbf{B},\mathbf{G})$ forms a $(1,2,\ldots,d-1)$-forbidding system.

Also, every $d$-set $S\in\cA$ is $(\mathbf{B},\mathbf{G})$-compatible, and $S^{(d)}$ consists of all the permutation of elements in $S$. Thus, $S^{(d)}_1,\ldots,S^{(d)}_m$ are mutually disjoint if $\cA=\{S_1,\ldots,S_m\}$. Therefore, We can apply \cref{theorem:entropykk}. Note that $|\cF|=d!|\cA|$ and $|\partial\cF|=(d-1)!|\partial\cA|$, and the statement follows.
\end{proof}

\begin{proof}[Proof of \cref{theorem:q-analog} using \cref{theorem:entropykk}]
Let $U=\mathbb{F}_q^n\setminus \{0\}$. Let $G_k$ be the family of all sets of $k$ linearly independent vectors. Let $B_k$ be the complement of $G_k$, i.e. $B_k$ consists of sets of $k$ linearly dependent vectors.

For any multiset $A$ with $k$ linearly dependent vectors, we know that any extension of $A$ is also linearly dependent. If $A\in G_k$, then there are exactly $q^{k}-1$ bad extensions of $A$. Namely, those extensions that extend by a vector in the span of $A$ are bad. Thus, the pair $(\mathbf{B},\mathbf{G})$ is a $(q-1,q^2-1,\ldots,q^{d-1}-1)$-forbidding system.
Also, for every $d$-dimensional subspace $S'\in\cA$, $S=S'\backslash\{0\}$ is $(\mathbf{B},\mathbf{G})$-compatible since it contains all $x$'s that are linearly dependent with $A$ if $\Span(A)$ is in $S'$.

If $\cA=\{S'_1,\ldots,S'_m\}$, we can define $\cF=\cup_i S_i^{(d)}$ where $S_i=S'_i\setminus\{0\}$. Note that any tuple in $S_i^{(d)}$ is a basis in $S'_i$ and thus uniquely determines $S'_i$, which shows that $S_1^{(d)},\ldots,S_m^{(d)}$ are mutually disjoint. Thus, We can apply \cref{theorem:entropykk}. Note that $|\cF|=(q^d-1)(q^d-q)\ldots(q^d-q^{d-1})|\cA|$ and $|\partial\cF|=(q^{d-1}-1)(q^{d-1}-q)\ldots(q^{d-1}-q^{d-2})|\partial\cA|$, and the statement follows.
\end{proof}

Our method is also capable of proving the multidimensional version of the Kruskal--Katona theorem considered by Bukh \cite{Bukh12}. For the actual statement, we refer the readers to that paper. However, the proof of this theorem using our entropy method is notationally heavy and most of the ideas are identical. We thus omit the proof.

\section{Other related problems inspired by joints}\label{section:OpenProblem}
In \cref{section:RainbowTriangle}, we considered the graph-theoretic counterpart of the multijoints problem and showed the best constant for $d=3$.
We can also consider the graph-theoretic counterparts of many other variants of the joints problem considered in \cite{TYZ22} and see if one could show a better constant in the graph-theoretic case.
Here, instead of stating the most general problem, we state some special cases where we still do not know the best possible constants.

For example, the graph-theoretic problem for joints of $2$-flats in $\RR^6$ is the following extremal problem regarding $4$-uniform hypergraphs.

\begin{question}
    Let $H$ be a $4$-uniform hypergraph with $N$ hyperedges, and let $J$ be the number of $6$-subsets $\Delta$ of $V(H)$ such that there are $\ell_1,\ell_2,\ell_3\in E(H)$ with $\Delta = (\Delta\backslash\ell_1)\sqcup (\Delta\backslash \ell_2)\sqcup(\Delta\backslash \ell_3)$.
    How large can $J^2/N^3$ be?
\end{question}
We emphasize that each $6$-subset is counted at most once in this setting.
The bound in \cite{TYZ22} translates directly into an upper bound of $10/3$ for this problem.
On the other hand, let $H$ be the $4$-uniform hypergraph on $\{x_1,x_2,x_3,x_4,y_1,y_2,y_3,y_4\}$ with the edges $\{x_1,x_2,x_3,x_4\}, \{y_1,y_2,y_3,y_4\}$ and $\{x_i,x_j,y_k,y_{\ell}\}$ for any $\{i,j\},\{k,\ell\}\subseteq \{1,2,3,4\}$ that are either equal or disjoint.
It is easy to verify that all $6$-subsets of $V(H)$ satisfy the condition and are thus counted in $J$, showing that $N=14$ and $J=28$.
This gives $J^2/N^3=2/7$, and we suspect that it is optimal.

The joints problems for flats of different dimensions translate into graph-theoretic problems where the hypergraph has mixed uniformity. Here is the graph-theoretic version of the joints problem where a joint is defined to be the intersection of a $2$-flat and two lines.

\begin{question}
    Let $H$ be a hypergraph with $N_2$ $2$-edges and $N_3$ $3$-hyperedges. We say that a $4$-subset $\Delta$ of $V(H)$ is \emph{good} if there exists a way to label $\Delta=\{v_1,v_2,v_3,v_4\}$ such that $\{v_1,v_2,v_3\},\{v_1,v_2,v_4\}$ are $3$-hyperedges and $\{v_3,v_4\}$ is a $2$-edge. Let $J$ be the number of good $4$-subsets. How large can $J^2/N_2N_3^2$ be?
\end{question}

A direct application of Shearer's inequality (\cref{theorem:ShearerStrong}) shows that the answer is at most $9/2$. 
To see this, uniformly sample $(X_1,X_2,X_3,X_4)$ from all quadruples of vertices such that $\{X_1,X_2,X_3\},\{X_1,X_2,X_4\},\{X_3,X_4\}$ are edges in $H$.
Note that for each good set, there are at least four labelings of the vertices that are sampled in this random process.
Thus, by Shearer's inequality,
\begin{align*}
    2\log_2(4J)&\leq 2H(X_1,X_2,X_3,X_4)\leq H(X_1,X_2,X_3)+H(X_1,X_2,X_4)+H(X_3,X_4)\\
    &\leq 2\log_2(6N_3)+\log_2(2N_2).
\end{align*}
This gives $J^2\leq 9N_2N_3^2/2$.

The bound in \cite{TYZ22} translates into an upper bound of $3$. 
There is a lower bound $3/2$ that is given by the following construction. 
Take three sets of vertices of size $n$ each. The $2$-edges are the edges that connect two vertices in the same part, and the $3$-hyperedges are the hyperedges that connect three vertices from all three parts. 
In this construction, we get $N_2=(1-o(1))3n^2/2, N_3=(1-o(1)) n^3$ and $J=(1-o(1)) 3n^4/2$ since all the $4$-sets that have at least a vertex in each part are good. Therefore, the ratio is $3/2$ in this construction. We again suspect that this construction is optimal.

To state the next question, we will consider a generalization of generically induced configurations.
In \cite{CY23}, we define \emph{projected generically induced configurations} as the image of some generically induced configurations in a higher-dimensional space via a general projection.
We may thus also consider the graph-theoretic version of the multijoints problem when we restrict to projected generically induced configurations.
In the case of $d=3$, we get the following question (see \cite{CY23} for details).

\begin{question}
    Let $H$ be a simple $(\delta+2)$-uniform hypergraph for some $\delta\in \ZZ_{\geq 0}$ with $R$ red edges, $B$ blue edges and $G$ green edges.
    Recall that a hypergraph is simple if no edges appear twice.
    Let $J$ be the number of $(\delta+3)$-subsets $\Delta$ of $V(H)$ such that it contains an edge of each color.
    How large can $J^2/RGB$ be?
\end{question}

We have already shown that the maximum of $J^2/RGB$ is $2$ if $\delta=0$.
However, for larger $\delta$, our argument no longer works.
Using results from joints (\cite[Theorem 3.2]{YZ23} and \cite[Section 7]{CY23}), we see that $J^2/RGB\leq 6$ for any $\delta$.
It seems that the stronger bound $J^2/RGB\leq 2$ could still hold for arbitrary $\delta$.

Lastly, we consider some purely graph-theoretic problems that can be resolved by the polynomial method, while there are still no known graph-theoretic proofs.
The first is the partial shadow problem proposed by Bollob\'as and Eccles \cite{BE15}.
To be consistent with the other problems stated here, we will state the problem in terms of hypergraphs.
Given positive integers $r,m$ and a nonnegative integer $k\leq r$, let $H$ be an $(r-1)$-uniform hypergraph.
The quantity $f(r,m,k)$ is the minimum number of edges in $H$ so that there are $m$ $r$-subsets of $V(H)$ that contain at least $r-k$ edges in $H$.
Using the connection between the partial shadow problem and the joints problem via projected generically induced configurations, we \cite{CY23} prove the following statement that was conjectured by Bollob\'as and Eccles \cite{BE15}.

\begin{theorem}
    Let $r,m$ be two positive integers, $k\leq r$ be a nonnegative integer and $x\geq r-k$ be a real number such that $m=\binom{x}{r-k}$.
    Then $f(r,m,k)\geq \binom{x}{r-k-1}.$
\end{theorem}

Note that when $k=0$, this is exactly Lovasz's version of the Kruskal--Katona theorem.
As noted in the introduction, a simple application of the Shearer's inequality immediately gives a bound that matches Lovasz's bound (and also the Kruskal--Katona bound) up to a $o(1)$-factor.
However, for $k$ strictly between $0$ and $r$, this is no longer the case.
To see this, we fix an $(r-1)$-uniform hypergraph $H$ with $m$ $r$-subsets of $V(H)$ that contain at least $r-k$ edges in $H$, and for each such $r$-subset $\Delta$, we pick exactly $r-k$ edges from the edges it contains.
Choose $\Delta$ uniformly at random, and choose the labeling $X_1,\ldots, X_r$ of the vertices in $\Delta$ uniformly so that the $r-k$ edges are $\Delta\backslash X_i$ for $i=k+1,\ldots, r$.
By Shearer's inequality,
\[(r-k-1)H(X_1,\ldots, X_r)\leq \sum_{i=k+1}^{r}H\left((X_j)_{j\neq i}\right).\]
The left hand side is exactly $(r-k-1)\log_2(k!(r-k)!\cdot m)$ as there are exactly $k!(r-k)!$ labelings for each $\Delta$.
The right hand side can be bounded using the uniform bound by $(r-k)\log_2((r-1)!\cdot e(H))$ as there are $(r-1)!$ different labelings for each edge.
Therefore
\[e(H)\geq \frac{\left(k!(r-k)!\cdot m\right)^{(r-k-1)/(r-k)}}{(r-1)!}=\left(\frac{k!^{(r-k-1)/(r-k)}}{(r-1)!}+o(1)\right)x^{r-k-1},\]
and we get a lower bound $f(r,m,k)\geq (k!^{(r-k-1)/(r-k)}/(r-1)! +o(1))x^{r-k-1}$ which has a worse constant than the bound $f(r,m,k)\geq (1/(r-k-1)!+o(1))x^{r-k-1}$ we get using the results from joints.
The reason why Shearer's inequality does not give a good bound here is that we are overcounting the number of labelings of the edges in the argument---we should only expect $(r-k-1)!$ different labelings for the tight example.

The next problem comes from translating the joints with multiplicity problem.

\begin{theorem}\label{theorem:MultiplicityGraphTheoretic}
    Let $d$ be a positive integer, and let $H$ be a complete weighted $(d-1)$-uniform hypergraph where the weights are nonnegative integers.
    Let $N$ be the total weight of the hyperedges.
    For any $d$-subset $\Delta$ of $V(H)$, define its weight $w(\Delta)$ to be the geometric mean of the weights of the $d$ hyperedges it contains.
    Then 
    \[\sum_{\Delta\in \binom{V(H)}{d}}w(\Delta)^{d/(d-1)} \leq \frac{\left((d-1)!\right)^{1/(d-1)}}{d}\cdot N^{d/(d-1)}.\]
\end{theorem}
Previously in \cite{TYZ22}, a weaker bound with the constant removed is proved (see (6) in the discussion after Theorem 1.10).
See \cite{CY23} for a full proof of this stronger result.
We include two other proofs here. 
The first one only works for $d=3$ by using spectral graph theory and the second one works for all $d$ by applying a powerful inequality from analysis.
\begin{proof}[Direct proof for $d=3$]
    In this case, the hypergraph $H$ is a weighted graph.
    Let $A$ be the matrix representing the weighted graph, where $a_{ij}$ is the weight of the edge connecting $i$ and $j$.
    Let $M$ be the term-wise square root of $A$.
    Then we can compute
    \[\tr(M^2) = 2N\]
    and
    \[\tr(M^3) = 6\sum_{\Delta\in\binom{V(H)}{d}}w(\Delta)^{3/2}.\]
    Let $\lambda_1,\ldots, \lambda_n$ be the eigenvalues of $M$.
    Then
    \[\tr(M^2)^3-\tr(M^3)^2 = \left(\sum_i\lambda_i^2\right)^3-\left(\sum_i\lambda_i^3\right)^2\geq \sum_{i<j}3\lambda_i^4\lambda_j^2-2\lambda_i^3\lambda_j^3+3\lambda_i^2\lambda_j^4\geq 0.\]
    Thus
    \[\sum_{\Delta\in\binom{V(H)}{d}}w(\Delta)^{3/2}\leq \frac{\sqrt{8}}{6}N^{3/2}=\frac{\sqrt{2}}{3}N^{3/2},\]
    as desired.
\end{proof}

Asaf Cohen Antonir kindly pointed out to us that it is possible to prove \cref{theorem:MultiplicityGraphTheoretic} using the so-called generalized H\"older's inequality \cite[Theorem 2.1]{Fin92}.
The version of the generalized H\"older's inequality we will use is the following.

\begin{theorem}[Generalized H\"{o}lder's inequality, discrete and simplified version]\label{theorem:GeneralizedHolder}
    Let $k$ be a positive integer.
    Let $I_1,\ldots,I_m$ be subsets of $[d]$ so that for each $i\in [d]$, $i\in I_j$ holds for at least $k$ different $j$'s. Let $f_i:\ZZ^{I_i}\rightarrow\mathbb{R}_{\geq 0}$ be a function for all $i\in [m]$. It follows that
    \[\sum_{x\in \ZZ^{d}}\prod_{i\in [m]}f_i(\pi_{I_i}(x))^{1/k}\leq \prod_{i\in [m]}\left(\sum_{y_i\in \ZZ^{I_i}}f_i(y_i)\right)^{1/k},\]
    where $\pi_{I_i}$ is the projection from $\ZZ^{d}$ to $\ZZ^{I_i}$.
    \end{theorem}
    \begin{remark}
        Generalized H\"{o}lder's inequality is a common generalization of H\"{o}lder's inequality and the Loomis--Whitney inequality, and it is closely related to the Brascamp--Lieb inequalities. The version we will use is actually the weighted version of the Loomis--Whitney inequality.
    \end{remark}
\begin{proof}[Proof of \cref{theorem:MultiplicityGraphTheoretic} using \cref{theorem:GeneralizedHolder}]
    Assume the vertex set of $H$ is $[n]$.
    We may apply \cref{theorem:GeneralizedHolder} by taking $I_i=[d]\setminus\{i\}$ for $i\in [d]$, and $k=d-1$. For all $i\in [d]$, let $f_i(v_1,\dots,\widehat{v_i},\dots,v_{d})$ be the weight of the hyperedge $\{v_1,\dots,\widehat{v_i},\dots,v_{d}\}$ and let $f_i$ be zero elsewhere. It follows that 
    \[\sum_{\Delta\in \binom{V(H)}{d}}d!\cdot w(\Delta)^{d/(d-1)} \leq \left((d-1)!\cdot N\right)^{d/(d-1)}\]
    since each $\Delta$ is counted $d!$ times and each hyperedge is counted $(d-1)!$ times over all possible permutations of the vertices. The theorem follows after rearranging.
\end{proof}

A possible approach of giving yet another proof of \cref{theorem:MultiplicityGraphTheoretic} is to translate the polynomial-method proof into a purely graph-theoretic one.
In general, it would be interesting to reprove the inequalities obtained from the vanishing lemmas (for example, Lemma 2.3 in \cite{YZ23}) in the purely graph-theoretic setting.
We believe that doing so would shed more lights on the problems considered here.

\section*{Acknowledgement}
    The authors would like to thank Noga Alon and Boris Bukh for useful comments on the paper, and Yufei Zhao for comments on the title.
    The authors would also like to thank Asaf Cohen Antonir for inspiring discussion that leads to the very last part of the paper.
    Lastly, the authors would like to thank the anonymous referees for their valuable suggestions.
\bibliographystyle{amsplain0}
\bibliography{ref_KK}

\appendix
\section{A construction with many rainbow tetrahedra}\label{appendix:tetrahedra}
The following is an explicit construction with $16$ red $3$-edges, $16$ blue $3$-edges, $12$ green $3$-edges, $12$ yellow $3$ edges and $48$ rainbow tetrahedra (i.e. $K^{(3)}_4$).
This would give a constant of
\[\frac{48^3}{16^2\cdot 12^2}=3.\]

To begin the construction, consider a vertex set $\{u_1,u_2,u_3,u_4,v_1,v_2,v_3,v_4\}$.
We list out the edges of each color below, where $u_5=u_1$ and $v_5=v_1$.
\begin{itemize}[leftmargin=*]
    \item Red: $\{u_i,u_{i+1},v_j\}$ for any $i,j\in [4]$ and $i\equiv j\mod 2$, and $\{u_i, v_j, v_{j+2}\}$ for any $i\in[4]$, $j\in [2]$.
    \item Blue: $\{u_i, v_i, v_{i+1}\}$ for any $i,j\in [4]$ and $i\equiv j\mod 2$, and $\{u_i, u_{i+2}, v_j\}$ for any $i\in[2]$, $j\in[4]$.
    \item Green: $\{u_i,u_{i+1},v_j\}$ for any $i,j\in[4]$ and $i\not\equiv j\mod 2$, and all the $3$-subsets of $\{v_1,v_2,v_3,v_4\}$.
    \item Yellow: $\{u_i,v_j,v_{j+1}\}$ for any $i,j\in[4]$ and $i\not\equiv j\mod 2$, and all the $3$-subsets of $\{u_1,u_2,u_3,u_4\}$.
\end{itemize}
It is easy to verify that each color gets the claimed number of $3$-edges, and the rainbow tetrahedra are $S_u\cup\{v_i\}$ for any $3$-subsets $S_u$ of $\{u_1,u_2,u_3, u_4\}$ and $i\in [4]$, $\{u_i\}\cup S_v$ for any $i\in[4]$ and $3$-subsets $S_v$ of $\{v_1,v_2,v_3, v_4\}$, and also $\{u_i,u_{i+1},v_j,v_{j+1}\}$ for any $i,j\in[4]$.
There are $4\cdot 4+4\cdot 4+4\cdot 4 = 48$ of them.
\end{document}